\documentclass[11pt]{amsart}
\usepackage{amssymb, graphics}
\usepackage{hyperref}
\newtheorem{lemma}{Lemma}
\newtheorem{prop}{Proposition}
\newtheorem{thm}{Theorem}
\newtheorem{cor}{Corollary}
\theoremstyle{definition}
\newtheorem{rem}{Remark}

\newcounter{numl}

\newcommand{\labelnuml}{\textup{(\roman{numl})}}
\newenvironment{numlist}{\begin{list}{\labelnuml}%
{\usecounter{numl}\setlength{\leftmargin}{0pt}%
\setlength{\itemindent}{2\parindent}%
\setlength{\itemsep}{\smallskipamount}\def
\makelabel ##1{\hss \llap {\upshape ##1}}}}{\end{list}}

\newcommand{\R}{{\mathbb R}}
\newcommand{\C}{{\mathbb C}}

\newcommand{\Z}{{\mathbb Z}}
\newcommand{\N}{{\mathbb N}}

\newcommand{\cL}{{\mathcal L}}

\newcommand{\cK}{{\mathcal K}}
\newcommand{\cO}{{\mathcal O}}

\newcommand{\symprod}{\mathbin{\raise1pt\hbox{$\scriptstyle\bigcirc$}}}

\newcommand{\delbar}{\overline{\partial}}
\newcommand{\del}{\partial}

\begin{document}

\title[Bihermitian metrics on Hopf surfaces]
{Bihermitian metrics on Hopf surfaces}

\author[V. Apostolov]{Vestislav Apostolov}
\address{Vestislav Apostolov \\ D{\'e}partement de Math{\'e}matiques\\
UQAM\\ C.P. 8888 \\ Succ. Centre-ville \\ Montr{\'e}al (Qu{\'e}bec) \\
H3C 3P8 \\ Canada \\ and \\
Institute of Mathematics and Informatics \\ Bulgarian Academy of Science \\ Acad. G. Bonchev Str. Bl. 8 \\ 1113 Sofia \\ Bulgaria}
\email{apostolov.vestislav@uqam.ca}
\author[G. Dloussky]{Georges Dloussky}
\address{Georges Dloussky \\ Centre de Math\'ematiques et d'Informatique \\ Laboratoire d'Analyse Topologie et probabilit\'es\\ Universit\'e d'Aix-Marseille 1\\ 39, rue Joliot-Curie, 13453 \\ Marseille Cedex 13 \\ France}
\email{dloussky@cmi.univ-mrs.fr}

\thanks{This paper was written while the first named author was visiting the University of Aix-Marseille~1. He is grateful to this institution for providing excellent working environment. Both authors wish to thank A.~Teleman for helpful discussions and for explaining his results in \cite{teleman}.}

\begin{abstract}
Inspired by a construction due to Hitchin~\cite{hitchin4}, we produce strongly bihermitian metrics on certain Hopf complex surfaces, which integrate the locally conformally K\"ahler metrics found by Gauduchon--Ornea~\cite{GO}. We  also show that the Inoue complex surfaces with $b_2=0$ do not admit bihermitian metrics. This completes the classification of the compact complex surfaces admitting strongly bihermitian metrics.  
\end{abstract}

\maketitle

\section{Introduction}

A {\it bihermitian structure} on a $4$-dimensional connected manifold $M$ consists of a pair of integrable complex structures $J_+$ and $J_-$ inducing the same orientation,  and a Riemannian metric $g$ which is hermitian with respect to both $J_{\pm}$. As a trivial example one can take a genuine hermitian structure $(g,J)$ and put $J_{\pm}=\pm J$. We exclude this situation by assuming that $J_+(x) \neq \pm J_-(x)$ at at least one point  $x$ of $M$. The special case when $J_+ \neq J_-$ {\it everywhere} on $M$ will be referred to as {\it strongly bihermitian} structure. 

In pure mathematics literature the theory of bihermitian $4$-manifolds was initiated~\cite{salamon,pontecorvo,kobak,AG0}  from the point of view of $4$-dimensional conformal geometry, where the two integrable almost-complex structures are given by the roots of the conformal Weyl tensor. Subsequent work~\cite{pontecorvo,AGG,A,dloussky} was mainly focussed on answering the following 

\vspace{0.2cm} 
\noindent {\it Question 1.} {When does a compact
complex surface $(M,J)$ admit a bihermitian structure $(g,J_+,J_-)$ with $J=J_+$?}

\vspace{0.2cm} There has been a great deal of interest in bihermitian geometry more recently, motivated by its link with the notion of {\it generalized K\"ahler geometry}, introduced
and studied by Gualtieri~\cite{gualtieri} in the context
of the theory of generalized geometric structures initiated by
Hitchin~\cite{hitchin}. It turns out~\cite{gualtieri} that a generalized K\"ahler structure is equivalent to the data of a Riemannian metric $g$ and two $g$-orthogonal complex structures  $(J_+,J_-)$, satisfying the relations
$d^c_+ F_+^g +d^c_- F_-^g = 0, \ dd^c_\pm F_\pm^g = 0,$ where $F^g_{\pm}(\cdot, \cdot) = g(J_{\pm}\cdot, \cdot)$ are the fundamental 2-forms of the hermitian structures $(g,J_{\pm})$, and $d^c_\pm$ are the associated
$i(\delbar_{\pm}-\del_{\pm})$ operators.   (These conditions on a pair of hermitian structures were, in fact, first described in the physics paper~\cite{physicists} as the
general target space geometry for a $(2,2)$ supersymmetric sigma
model.) In four dimensions we obtain a bihermitian structure, provided that the generalized K\"ahler structure is of {\it even} type (which corresponds to our assumption that $J_{\pm}$ induce the same orientation on $M$) and $F_+^g\neq \pm F_-^g$ (i.e. $J_+ \neq \pm  J_-)$. 

When the first Betti number of $M$ is {\it even}, it follows from \cite{AGG} (see also Remark~\ref{gk-geometry} below) that,  conversely,  any bihermitian metric arises from a generalized K\"ahler structure, up to a conformal change of the metric. This correspondence was used to produce a number of new constructions of bihermitian metrics~\cite{goto,gualtieri-et-al,hitchin2,Lin-Tolman}, where the existence of a K\"ahler structure on $M$ plays a crucial r\^ole.  

\vspace{0.2cm} When the first Betti number $b_1(M)$ is {\it odd}, on the other hand, the manifold $(M,J)$ does not admit K\"ahler metrics, nor a bihermitian conformal structure necessarily  come from a generalized K\"ahler structure. (As a matter of fact no strongly bihermitian metric  can be obtained from a generalized K\"ahler structure, see Remark~\ref{gk-geometry}). Thus, examples of bihermitian structures in this case are still scarce (see, however, \cite{pontecorvo,AGG}).

\vspace{0.2cm}Thus motivated, we prove in this note the following result.
\begin{thm}~\label{main} A compact complex surface $(M,J)$ with odd first Betti number admits a strongly bihermitian structure $(g,J_+,J_-)$ with $J_+=J$ if and only if it is a Hopf surface whose canonical bundle defines a class in the image of  $H^1(M, \R_+^*) \to H^1(M, \cO^*)$.

These are the Hopf surfaces with universal covering space $\C^2 \setminus
\{(0,0)\}$ and fundamental group $\Gamma$  which belongs to one of the following cases:
\begin{enumerate}
\item[\rm (a)] $\Gamma$ is generated by the automorphism $(z_1,z_2) \mapsto (\alpha z_1, {\bar \alpha} z_2)$ {\rm (}where $\alpha \in \C, 0< |\alpha|<1${\rm )}, and a  finite subgroup of $SU(2)$;
\item[\rm (b)] $\Gamma$ is generated by the  automorphisms $(z_1,z_2) \mapsto (\alpha z_1, a\alpha^{-1} z_2)$ {\rm (}where $\alpha \in \C,  0<|\alpha|^2<a<|\alpha| <1${\rm )}, and   $(z_1,z_2)
\mapsto (\varepsilon z_1, \varepsilon^{-1} z_2)$ with $\varepsilon$ primitive
$\ell$-th root of $1$;
\item[\rm (c)] $\Gamma$ is generated by the automorphisms 
 $(z_1,z_2) \mapsto (\beta^m z_1 + \lambda z_2^m, \beta  z_2)$ {\rm (}where
$\lambda \in {\mathbb C}^*, \ \beta^{m+1} = a,  \ 0<a <1${\rm )}, and $(z_1,z_2) \mapsto (\varepsilon z_1, \varepsilon^{-1} z_2)$ with $\varepsilon$ primitive $\ell$-th root of $1$, where $m = k\ell -1$,  $k\in {\mathbb N}^*$.
\end{enumerate}
\end{thm}

Note that the Hopf surfaces with fundamental group belonging to the case (a) of the above theorem are precisely those which admit a compatible hyperhermitian structure (see~\cite{boyer}). A description of all possible finite subgroups of $SU(2)$ that appear in this case can be found in \cite{kato1}.

Combined  with the results in \cite{AGG}, Theorem~\ref{main} yields the complete list of the compact complex surfaces $(M,J)$ which carry  a strongly bihermitian structure with $J_+=J$: these are $K3$ surfaces, complex tori,  and the Hopf complex surfaces of the type described above. Note that in the case when $b_1(M)$ is even, $(M,J)$ does also admit a compatible hyperhermitian structure~\cite{boyer}, but this is not longer true when $b_1(M)$ is odd.

\vspace{0.2cm}Compared to the partial results in \cite{A}, Theorem~\ref{main} brings in two new ingredients. 

First of all, we show  that the Inoue surfaces with $b_2(M)=0$ (defined and studied in \cite{inoue}) do not admit compatible bihermitian structures. This is achieved by combining a recent  result of A.~Teleman~\cite[Rem.~4.2]{teleman} (according to which on any  Inoue surface the degree of the anti-canonical bundle is negative with respect to any standard hermitian metric) together with an observation from \cite[Rem.~1]{A} (that this degree must be positive, should a bihermitian structure exists).

Secondly, we address the question of existence of bihermitian structures on the Hopf surfaces in our list. We give a concrete construction of bihermitian metrics by using an idea of Hitchin~\cite{hitchin4}. It allows us to deform the locally conformally K\"ahler metrics $(g_{GO},J)$ found by Gauduchon--Ornea~\cite{GO} in order to obtain a family of bihermitian metrics $(g_t, J, J_-^t)$  with $J_-^t= \phi_t^*(J)$, where $\phi_t$ is a path of diffeomorphisms;  as $t\to 0$, $J_-^t \to J$ and $g_t/t \to g_{GO}$.

\section{Bihermitian geometry -- preliminary results}
In this section we recall some properties of bihermitian
metrics, which we will need for the proof of our main results. The articles \cite{A,AGG,hitchin4} are  relevant references for more details. 

Let $(g,J_+, J_-)$ be a bihermitian structure on a $4$-manifold $M$, i.e.
a riemannian metric $g$ and two $g$-compatible complex
structures $J_\pm$ with $J_+(x)\neq \pm J_-(x)$ at some point $x\in M$, and such that  $J_+$ and $J_-$ induce the same orientation on $M$. Notice that any riemannian metric which is conformal to $g$ is again bihermitian with respect to $J_{\pm}$; we can therefore define a bihermitian {\it conformal} structure $(c=[g],J_+,J_-)$ on $M$. 

For a fixed metric $g \in c$, we consider the fundamental $2$-forms $F_+^g$ and $F_-^g$ of the hermitian metrics $(g,J_+)$ and $(g,J_-)$,  defined by 
$F_{\pm}^g(\cdot, \cdot)=g(J_{\pm}\cdot, \cdot).$
The corresponding Lee $1$-forms, $\theta^g_{+}$ and $\theta_-^g$, are  introduced by the relations $dF_+^g = \theta_+^g\wedge F_+^g$ and $dF_-^g= \theta_-^g\wedge F_-^g$, or equivalently $\theta_{\pm}^g = J_{\pm} (\delta^g F_{\pm}^g)$,  where  $\delta^g$ is the co-differential of $g$ and  the action of an almost complex structure $J$ on  the cotangent bundle $T^*(M)$ is given by $(J\alpha)(X) = - \alpha(JX)$. Notice that with respect to a conformal change of the metric $\tilde g= e^f g$, the Lee forms change by  $\theta_{\pm}^{\tilde g} = \theta_{\pm}^g + df$. 

The following result was established in \cite{AG0,AGG}.
\begin{lemma}\label{lee-forms} For any bihermitian metric $(g,J_+,J_-)$, the Lee forms $\theta_{\pm}^g$ satisfy 
\begin{equation*}
2\delta^g \theta_+^g + |\theta_+^g|_g^2 = 2\delta^g \theta_-^g + |\theta_-^g|_g^2,
\end{equation*}
where  $|\cdot |_g$ is the point-wise norm induced by $g$.

Moreover, the $2$-form $d(\theta^g_+ + \theta^g_ -)$ is anti-selfdual so that, when $M$ is compact, $\theta^g_+ + \theta^g_-$ is closed.
\end{lemma}

The complex structures $J_{\pm}$ must satisfy the relation \cite{pontecorvo,AGG}
\begin{equation}\label{anticommute}
J_+ J_- + J_- J_+ = -2p  \ {\rm Id},
\end{equation}
where $p=-\frac{1}{4}{\rm trace}(J_+J_-)$ is the so-called {\it angle function} which, at each point $x \in M$,   verifies $|p(x)|\le 1$ with $p(x)= \pm 1$ if and only if $J_+(x) = \pm J_-(x)$. In particular, in the strongly bihermitian case, $|p|<1$ everywhere.

Another natural object associated to the bihermitian structure $(g,J_+,J_-)$ is the $g$-skew-symmetric   endomorphism $[J_+,J_-]=J_+J_--J_-J_+$ which anti-commutes with both $J_{\pm}$. Using the metric $g$, we define a real $J_{\pm}$-anti-invariant $2$-form by
$$\Phi^g(\cdot, \cdot)= \frac{1}{2}g ([J_+,J_-] \cdot, \cdot).$$
Letting  
$$\Psi_{\pm}^g(\cdot, \cdot):=-\Phi^g(J_{\pm} \cdot, \cdot),$$
we consider  the complex $2$-forms
$$\Omega_{\pm}^g(\cdot, \cdot) = \Phi^g + i \Psi_{\pm}^g, $$
which are of type $(2,0)$ on the respective complex manifolds $(M,J_{\pm})$. 

Using \eqref{anticommute} and that $F^g_{\pm}, \Phi^g, \Psi_{\pm}^g$ are all selfdual $2$-forms satisfying the obvious orthogonality relations given by their types with respect to $J_{\pm}$, it is straightforward to check (see \cite[Lemma~2]{AGG})
\begin{equation}\label{two-forms}
\begin{split}
& F_+^g = p F_-^g + \Psi^g_-, \  \  \ F_-^g = pF_+^g - \Psi_+^g,  \\ 
& {\Phi^g}\wedge \Phi^g  = \Psi^g_{\pm}\wedge \Psi_{\pm}^g= 2(1-p^2)dv_g,  \ \    \Phi^g\wedge \Psi_{\pm}^g=0,  \ \ \Psi_+^g \wedge \Psi_-^g = p \Phi^g\wedge \Phi^g, \\
& (\Psi_-^g)^{1,1}_{J_+}=(1-p^2) F_+^g,
\end{split}
\end{equation}
where $(\cdot )^{1,1}_{J_+}$ denotes the $J_+$-invariant part of a two form.

Still using the metric $g$, we define the $g$-duals of the $(0,2)$-forms ${\bar \Omega}_{\pm}^g$, which are smooth sections, say $\sigma_{\pm}^g$, of the respective anti-canonical bundles ${\mathcal K}^{-1}_{J_{\pm}}=\wedge^2(T^{1,0}_{J_{\pm}}(M))$. Notice that at any point where $J_+(x) \neq \pm J_-(x)$,  $\Omega_{\pm}^g$ and $\sigma_{\pm}^g$ do not vanish, because of \eqref{anticommute}.

The following key result was established in \cite[Lemmas 2,3]{AGG}.
\begin{lemma}\label{poisson-structure} For any bihermitian metric $(g,J_+,J_-)$ the $(2,0)$-forms $\Omega_{\pm}^g$ and bi-vector fields $\sigma_{\pm}^g$ satisfy  
\begin{equation*}
\begin{split}
d\Omega_{\pm}^g &= \Big(\frac{1}{2}(\theta^g_+ + \theta^g_-) + d \log (1-p^2)\Big) \wedge \Omega_{\pm}^g,  \\
{\bar \partial}_{\pm} \sigma^g_{\pm} &=  -\frac{1}{2}(\theta_+^g + \theta_-^g)_{J_\pm}^{0,1}\otimes \sigma^g_{\pm}, 
\end{split}
\end{equation*}
where the first equality holds on the open subset of points where $J_+(x) \neq \pm J_-(x)$, and in the second equality ${\bar \partial}_{\pm}$ stand for the corresponding Cauchy--Riemann operators on the respective anti-canonical bundles $\cK^{-1}_{J_{\pm}}$, and $(\cdot )^{0,1}_{J_\pm}$ denote the $(0,1)$-parts taken relatively to the respective complex structure.
\end{lemma}

The first equality in the above lemma tells us that the real $2$-forms $\Phi^g, \Psi_{\pm}^g$ satisfy 
$d\Phi^g = \tau \wedge \Phi^g, d\Psi_{\pm}^g = \tau \wedge \Psi_{\pm}^g$ with $\tau=\frac{1}{2}(\theta^g_+ + \theta^g_-) + d \log (1-p^2)$. Together with the relations \eqref{two-forms}, this allows to reconstruct the strongly bihermitian metric from the $2$-forms  $\Phi^g, \Psi_+^g, \Psi_-^g$.

\begin{prop}\label{criterion}  Suppose $\Phi, \Psi_+, \Psi_-$ are non-degenerate real $2$-forms on $M$ which satisfy the  relations
\begin{equation*}
\begin{split}
\Phi^2 &= \Psi_+^2= \Psi_-^2, \ \  \Phi \wedge \Psi_{\pm} =0,  \ \  \\
d\Phi &= \tau \wedge \Phi, \ \ d\Psi_{\pm} = \tau \wedge \Psi_{\pm},
\end{split}
\end{equation*}
for some $1$-form $\tau$. 

Then, there exists a strongly bihermitian metric $(g,J_+,J_-)$ on $M$ with $\Phi^g=\Phi, \ \Psi_{\pm}^g=\Psi_{\pm}$ if and only if 
$\Psi_+\wedge \Psi_-= p \Phi^2,$
where $p$ is a smooth function with $|p|<1$
\end{prop}
\begin{proof} The necessity of the condition follows from \eqref{two-forms} and Lemma~\ref{poisson-structure}. (Recall that $p(x)= \pm 1$ if and only if $J_+(x) = \pm J_-(x))$.

In the other direction the result was originally established in \cite[Thm.~2]{AGG} in the case when $\tau=0$. Following \cite{AGG},  the almost complex structures $J_{\pm}$ are introduced by $\Psi_{\pm}(\cdot, \cdot)=-\Phi(J_{\pm} \cdot, \cdot)$,  while the conformal structure is determined by the property that $\Phi, \Psi_+, \Psi_-$ are selfdual $2$-forms.  The only notable difference with the proof given in \cite[Thm.~2]{AGG} is in establishing the integrability of $J_{\pm}$. For that we would need to generalize \cite[Lemma~5]{AGG} in order to show that for any pair  of $2$-forms, $(\Phi, \Psi)$, satisfying $$\Phi^2=\Psi^2, \ \Phi\wedge \Psi =0,$$ the almost-complex  structure they define by $\Psi(\cdot, \cdot)  = -\Phi(J\cdot, \cdot)$  is {\it integrable}, provided that $$d\Phi= \tau \wedge \Phi, \ \ d\Psi = \tau \wedge \Psi,$$
for some $1$-form $\tau$. This fact is well-known (see e.g. \cite{salamon}) but we give here a short proof for completeness.  Note that, by the assumption made, the differential of the $(2,0)$-form $\Omega=\Phi + i \Psi$ does not have a component of type $(1,2)$. Thus,  for any two complex vector fields  ${\bar V}, {\bar W}$ of type $(0,1)$, and any complex vector field $U$ of type $(1,0)$,  we have
$$0=d\Omega (U,{\bar V}, {\bar W})= -\Omega([{\bar V}, {\bar W}]^{1,0},U).$$
As $\Omega$ is a non-degenerate $2$-form on the complex vector bundle $T^{1,0}_J(M)$, it follows that  $[{\bar V}, {\bar W}]^{1,0}=0$, i.e. $J$ is integrable.\end{proof}

In the case when $J_+=J$ is given, there is a useful ramification of the above criterion, due to Hitchin.
\begin{cor}\label{hitchin-lemma}~\cite{hitchin4} Let $(M,J)$ be a complex surface endowed with a non-vanishing complex $(2,0)$-form $\Omega=\Phi + i \Psi$ such that $d\Omega = \tau \wedge \Omega$ for some real $1$-form $\tau$. 

Then, any real $2$-form $\Psi_-$ on $M$ satisfying $\Psi_-^2 = \Phi^2, \ \ d\Psi_- = \tau \wedge \Psi_-$, and whose $(1,1)$-part with respect to $J$ is positive definite, gives rise to a strongly bihermitian metric $(g,J_+,J_-)$ on $M$ with $J_+=J$.
\end{cor}
\begin{proof} Let $F: = (\Psi_-)^{1,1}_J$ be the positive definite $(1,1)$-component  of $\Psi_-$ with respect to the complex structure $J$. By putting $g(\cdot, \cdot) = F(\cdot, J\cdot)$, we obtain a Riemannian metric with respect to which $\Phi, \Psi_+=\Psi$ and $\Psi_-$ are selfdual $2$-forms of equal length. It then follows that $\Psi_+ \wedge \Psi_- = p \Phi \wedge \Phi$ for a smooth function $p$ satisfying $|p(x)|\le 1$ with $p(x)= \pm 1$ if and only if $\Psi_+(x) = \pm \Psi_-(x)$. The later inequality is impossible because the $(1,1)$-part of $\Psi_-$ is positive definite while $\Psi_+$ is $J$-anti-invariant. (This is also consistent with the last identity in \eqref{two-forms}.) We can now apply  Proposition~\ref{criterion}. \end{proof}

On a {\it compact} bihermitian $4$-manifold $(M,g,J_+,J_-)$, an overall assumption which we will make from now on, Lemma~\ref{poisson-structure} has  the following interpretation. 

To simplify notation,  consider one of the complex structures, say $J=J_+$, and drop the index $+$ for the corresponding $(2,0)$-form $\Omega^g$, bi-vector field $\sigma^g$, canonical bundle ${\mathcal K}_J$, Cauchy--Riemann operator $\bar \partial$,  etc. 

Since $\frac{1}{2}(\theta_+^g + \theta_-^g)$ is a closed $1$-form by Lemma~\ref{lee-forms}, it defines a holomorphic structure on the trivial (smooth) complex bundle $M\times \C$, by introducing the  connection $\nabla := \nabla^0 + \frac{i}{2}J(\theta_+^g + \theta_-^g)$, where $\nabla^0$ is the flat connection on $M\times \C$;   the new holomorphic structure depends, in fact,  only on the de Rham class $a=\frac{1}{2}[\theta_+^g + \theta_-^g]$, so we denote this holomorphic  line bundle by ${\mathcal L}_{a}$. This is consistent with the sequence of natural morphisms
\begin{equation}\label{morphisms} 
H^1(M,\R)\mapsto H^1(M, \R_+^*) \hookrightarrow  H^1(M, \C^*) \mapsto H_0^1(M,\cO^*),
\end{equation}
where the first morphism is induced by the exponential map (and hence is an isomorphism),  and $H_0^1(M,\cO^*)$ is the space of equivalent classes of topologically trivial holomorphic bundles. We will say that a topologically trivial holomorphic line bundle $\cL \in H_0^1(M,\cO^*)$  is of {\it real type} if it is in the image of $H^1(M, \R_+^*)$.

The second equality of Lemma~\ref{poisson-structure} simply means that $\sigma^g$ is a holomorphic section of  ${\mathcal K}^{-1}_J\otimes \cL_{-a}$, while the first equality and \eqref{two-forms} imply that its inverse, $\frac{1}{|\sigma^g|^2_g} {\Omega}^g$, is a meromorphic section of $\cK_J\otimes \cL_a$ (where $\cL_{-a}= \cL_{a}^{-1}$ in terms of the morphisms \eqref{morphisms}). Notice that by its very definition, $\sigma^g(x)=0$ if and only if $J_+(x)$ and $J_-(x)$ commute, which in view of \eqref{anticommute}, means $J_+(x)= \pm J_-(x)$. In other words,  ${\mathcal K}^{-1}_J\otimes \cL_{-a} \cong \cO$ if and only if $(c,J_+,J_-)$ is strongly bihermitian.

We now recall the definition~\cite{gauduchon0} of {\it degree} of a holomorphic line bundle with respect to a {\it standard} hermitian metric on $(M,J)$. A hermitian metric $g$ on a compact complex surface $(M,J)$ is called standard if its fundamental $2$-form $F^g$ is $dd^c$-closed (equivalently $\partial {\bar \partial} F^g=0$); this is the same as requiring that the corresponding Lee form $\theta^g$ is co-closed. A fundamental result of Gauduchon (see e.g.~\cite{gauduchon1}) affirms that such a metric exists, and is unique up to homothety, in each conformal class $c$ of hermitian metrics on $(M,J)$. Given a standard metric $g$ on $(M,J)$,  one defines the degree of a holomorphic line bundle $\cL$ by
\begin{equation}\label{degree}
{\rm deg}_g(\cL) = \frac{1}{2\pi}\int_M \rho \wedge F^g,
\end{equation}
where $\rho$ is the real curvature form of any holomorphic connection on $\cL$ (with the usual convention that  $\rho$ represents $2\pi c_1(\cL)_{\R}$).  Thus defined, it has the usual properties of  degree, notably it gives the volume with respect to $g$ of the divisor defined by any meromorphic section of $\cL$; in particular, the degree is non-negative for bundles with holomorphic sections, and is positive if there are sections with zeroes.

For a topologically trivial holomorphic line bundle of real type, say  $\cL_{a}$ ($a \in H^1_{dR}(M,\R)$), its degree with respect to a standard metric $g$ is easy to compute: we choose a representative closed $1$-form $\xi$ for $a$  and consider the holomorphic connection  $\nabla^{\xi} = \nabla^0 + iJ\xi$ on $\cL_a$, with curvature $\rho = d^c\xi$. Substituting in  \eqref{degree} and integrating by parts gives
\begin{equation}\label{harmonic-degree}
{\rm deg}_g (\cL_a)= \frac{1}{2\pi} \langle \xi, \theta^g \rangle_g=\frac{1}{2\pi}\langle a^g_h, \theta^g_h \rangle_g,
\end{equation}
where $a^g_h$ is the harmonic representatives of $a$ with respect of $g$, $\theta^g_h$ is the harmonic part of the Lee form $\theta^g$ and $\langle \cdot, \cdot \rangle_g$ is the global $L^2$ product on $1$-forms induced by the standard metric $g$. (For the last equality we have used that $\theta^g$ is co-closed.)

In the case of a bihermitian conformal structure $(c,J_+,J_-)$ on $(M,J)$ (with $J=J_+$), we take a standard metric  $g$ in $c$ with respect to $J$ and calculate the degree of $\cL_a$ with $a=\frac{1}{2}[\theta_+^g + \theta_-^g]$ as above.  By using  the first relation of Lemma~\ref{lee-forms} one gets (see \cite{A})
\begin{equation}\label{degree-computation}
{\rm deg}_g (\cL_{a})= \frac{1}{4\pi} \langle \theta_+^g + \theta_-^g, \theta_+^g \rangle_g= \frac{1}{8\pi}||\theta_+^g + \theta_-^g||_{L^2}^2.
\end{equation}

Notice that by \cite{gauduchon1},  the harmonic part of $\theta^g_h$ is zero if $b_1(M)$ is even. In this case the degree of $\cL_a$ must be zero,  hence $\cL_a = \cO$ and $\theta_+^g + \theta_-^g=0$ with respect to the standard metric of $(c,J_+)$ (which is therefore standard  too with respect to $J_-$), cf.~\cite[Lemma 4]{AGG}.

As $\sigma^g$ is a non-zero holomorphic section of ${\mathcal K}^{-1}_J\otimes \cL_{-a}$, this bundle has always positive or zero degree (the later is possible only when the bundle is trivial, i.e. the metric is strongly bihermitian). We thus obtain that ${\rm deg}_g ({\mathcal K}^{-1}_J)$ is either positive or is zero (and this is true for any standard metric on $(M,J)$). The zero case can only happen if both ${\mathcal K}_J$ and $\cL$ are trivial, and $\theta_+^g + \theta_-^g=0$ with respect to the standard metric $g$ of the bihermitian conformal structure. However, the later combination is impossible when $b_1(M)$ is odd, as shown in \cite[Prop.~4]{AGG}.

We summarize our discussion in the following  proposition which  gathers most of the information used so far to narrow the list of compact complex surfaces which can possibly admit bihermitian metrics \cite{A,AGG,dloussky,pontecorvo}.

\begin{prop}~\cite{A,AGG}\label{obstruction} Let $(M,J)$ be a compact complex surface. Any bihermitian conformal structure $(c,J_+,J_-)$ with $J=J_+$ defines a topologically trivial holomorphic line bundle ${\cL}$ of real type and non-positive degree with respect to the standard metric of $(c,J)$, and a non-zero holomorphic section of ${\mathcal K}^{-1}_J\otimes \cL$. In particular, the degree of ${\mathcal K}^{-1}_J$ must be non-negative with respect to the standard metric of $(c,J)$. 

Moreover, $\cL$ is trivial when $b_1(M)$ is even and the degree of ${\mathcal K}^{-1}_J$ is positive when $b_1(M)$ is odd.

The holomorphic bundle ${\mathcal K}^{-1}_J\otimes \cL$ is trivial if and only if $(g,J_+,J_-)$ is strongly bihermitian.
\end{prop}
\begin{rem}\label{sign} One can show that on a compact bihermitian surface the signs of  ${\rm deg}_g(\cL)$ and ${\rm deg}_g(\cK)$ are the same with respect to {\it any} standard metric on $(M,J)$. 

In the case when $b_1(M)$ is even this is obvious because, by Proposition~\ref{obstruction},  $\cL= \cO$ and $H^0(M, \cK^{-1}) \neq 0$.

In the case when $b_1(M)$ is odd, Proposition~\ref{obstruction} tells us that the degree of the anti-canonical bundle (with respect to a particular standard metric) is positive and, therefore, by  the properties of degree,  the pluricanonical line bundles ${\mathcal K}_J^{\otimes m}$ ($m \ge 1$) do not have non-zero holomorphic sections, i.e.~$(M,J)$ is of Kodaira dimension $-\infty$. This means that $(M,J)$ belongs to the class VII of the Enriques--Kodaira classification~\cite{bpv}. In particular, $b_1(M)=1$. According  to \cite{gauduchon2}~\footnote{The result used here is readily available in \cite[Prop.~1]{AG}.}, on a compact complex surface with odd first Betti number  the harmonic part of  the Lee form  $\theta_h^g$ of  {any} standard metric $g$ is non-zero. Since the space of standard metrics is convex (and $b_1(M)=1$), we obtain from \eqref{harmonic-degree} that ${\rm deg}_g (\cL_a)$ has the same sign (positive or zero by \eqref{degree-computation}) for {\it any} standard metric $g$ on $(M,J)$. The same is true for ${\rm deg}_g(\cK^{-1})$ because, by Proposition~\ref{obstruction}, $H^0(M, \cK^{-1}\otimes \cL)\neq 0$.\end{rem}

\begin{rem}~\label{gk-geometry} A bihermitian metric $(g,J_+,J_-)$ corresponds to a (twisted) generalized K\"ahler structure on $TM\oplus T^*M$ in the sense of \cite{gualtieri} if it verifies the extra relations
$$d^c_{+} F^g_{+} + d^c_-F^g_-=0,  \ \ dd^c_{\pm} F_{\pm}^g=0.$$
Since $J_{\pm}$ yield the same orientation on $M$, the first relation is equivalent (by applying the Hodge $*$ operator) to $\theta_+^g + \theta_-^g=0$,  while the second relations mean that $g$ is standard metric with respect to both $J_{\pm}$ (i.e. $\delta^g \theta_{\pm}^g=0$). Thus, for a bihermitian metric to come from a generalized K\"ahler structure the line bundle $\cL$ in Proposition~\ref{obstruction} must be trivial. Conversely, any bihermitian conformal structure $(c,J_+,J_-)$ with $\cL= \cO$ satisfies the above conditions with respect to the standard metric of $(c,J_+)$ (and thus, by \cite[Thm.~6.37]{gualtieri},  induces a twisted generalized K\"ahler structure on $TM \oplus T^*M$). This follows from the formula~\eqref{degree-computation} which shows that $\theta_+^g + \theta_-^g=0$ with respect to the standard metric of $(c,J_+)$, provided that ${\rm deg}_g (\cL)=0$.   Proposition~\ref{obstruction} implies that a strongly bihermitian surface with odd first Betti number is never of this type. \end{rem}

\section{Strongly bihermitian metrics on complex surfaces with odd first Betti number--Proof of Theorem~\ref{main}} We now specialize  to the case when $(M,J)$ has odd first Betti number and admits a strongly bihermitian metric. 

As we explained in Remark~\ref{sign}, the fact that $(M,J)$ has a standard metric with respect to which the degree of the anti-canonical bundle is positive (see Proposition~\ref{obstruction}) implies  that $(M,J)$ belongs to the class VII of the Enriques--Kodaira classification~\cite{bpv}. In particular,  $b_1(M)=1$ and $b^+(M)=0$,  where $b^{+}(M)=\frac{1}{2}(b_2(M) + \sigma(M))$ is the dimension of the space of self-dual harmonic $2$-forms with respect to any metric on $M$, and $\sigma(M)$ is the topological signature. In the strongly bihermitian case  ${\mathcal K}_J = \cL$ is topologically trivial (see Proposition~\ref{obstruction}), and then
$$0=c_1^2(M)= 2e(M) + 3\sigma(M)= -b_2(M),$$
where $e(M)= 2 - 2b_1(M) + b_2(M)$ is the Euler characteristic. Thus, $(M,J)$ is  a minimal surface of class VII with zero second Betti number. By Bogomolov's theorem (see~\cite{teleman0} or \cite{yau2} for a proof), $(M,J)$ is  either a minimal Hopf surface or an Inoue surface described in \cite{inoue}.

Combining Proposition~\ref{obstruction} with a result of \cite[Rem.~4.2]{teleman} allows us to exclude the case of Inoue surfaces. We reproduce the proof for completeness.

\begin{lemma}\cite{teleman} The degree of the anticanonical bundle of an Inoue complex surface with zero second Betti number is negative with respect to any standard metric. By Proposition~\ref{obstruction},  such a complex surface can not admit compatible bihermitian metrics.
\end{lemma}
\begin{proof}  
Following  \cite{teleman}, we show that the {\it canonical} bundle of an Inoue complex surface has positive degree with respect to any standard metric.

Recall that an Inoue surface is a quotient of ${\mathbb H} \times \C$ (where ${\mathbb H}$ denotes the upper half-plane) by properly discontinuous group $\Gamma$ of affine transformations. An explicit description of all such groups can be found in \cite{inoue}, and they belong to one the following three types.

For the Inoue surfaces $S_M$, the generators of $\Gamma$ are of the form
$$\gamma_0(w,z) = (\alpha w, \beta z), \ \gamma_k(w,z) = (w + a_k, z+b_k), \ k=1,2,3,$$
with $\alpha, a_k \in \R$ and $\alpha |\beta|^2 =1$ (where $w \in {\mathbb H}, z \in \C$).   The tensor 
$${\rm Im}(w) \Big(\frac{\partial}{\partial w} \wedge \frac{\partial}{\partial z}\Big) \otimes \Big(\frac{\partial}{\partial \bar w} \wedge  \frac{\partial}{\partial \bar z}\Big)$$
is  clearly invariant under $\Gamma$, and therefore defines  (after symmetrising) an hermitian metric on the canonical bundle of $S_M$. The curvature of this metric is $-i\partial {\bar \partial} \log ({\rm Im}(w))= \frac{i}{{\rm Im}(w)^2} dw\wedge d {\bar w}$,  which is a non-negative (but non-zero) $2$-form, so the degree of the canonical bundle is positive with respect to any standard metric.

Similarly, for the other two types of Inoue surfaces, $S^{+}_{N,p,q,r,t}$ and $S^-_{P,p,q,r}$, the generators of $\Gamma$ are of the form 
$$ \gamma_0(w,z)= (\alpha w, \epsilon z + t), \ \gamma_k(w,z)= (w+a_k, z+b_k w + c_k), \ k=1,2,3, $$
with $\alpha, a_k, b_k, c_k \in \R$ and $\epsilon=\pm 1$. Now a  $\Gamma$-invariant tensor is 
$${\rm Im}(w)^2 \Big(\frac{\partial}{\partial w} \wedge \frac{\partial}{\partial z}\Big) \otimes \Big(\frac{\partial}{\partial \bar w} \wedge  \frac{\partial}{\partial \bar z}\Big)$$
which defines an hermitian metric on the canonical bundle with  non-negative curvature $\frac{2i}{{\rm Im}(w)^2} dw\wedge d {\bar w}$; one concludes  as before that  the degree of the canonical bundle is positive with respect to any standard metric. \end{proof}

Let us now examine the Hopf surfaces. These are, by definition, compact complex surfaces with universal covering space $\C^2\setminus \{(0,0)\}$. It was shown by Kodaira~\cite{kodaira}  that  the fundamental group $\Gamma$ of such a surface is a finite extension of the infinite cyclic group $\Z$. The list of concrete realizations of $\Gamma$ as a group of automorphisms of $\C^2$ can be found in \cite{kato} and we shall make extensive use of this classification in the following rough form:  $\Gamma \cong \langle \gamma_0 \rangle \ltimes H$,  where 
\begin{enumerate}
\item[$\bullet$] $\langle \gamma_0 \rangle$ denotes the infinite cyclic group generated by a contraction of $\C^2$ of the form
\begin{equation}\label{contraction}
\gamma_0(z_1,z_2)= (\alpha z_1 + \lambda z_2^{m}, \beta z_2), 
\end{equation}
where $\alpha, \beta, \lambda \in \C$ with $ 0 <|\alpha| \le |\beta| <1, \ \lambda(\alpha - \beta^m)=0, m \in \N^*$; 
\item[$\bullet$] $H$ is  a finite subgroup of $U(2)$, subject to the following constraint
\begin{enumerate}
\item[\rm (i)] if $\lambda \neq 0$, then $H\subset U(1)\times U(1)$ is abelian and commutes with $\gamma_0$;
\item[\rm (ii)] if $\lambda=0$ and $|\alpha| \neq |\beta|$, then $H \subset U(1)\times U(1)$.
\end{enumerate}

\end{enumerate}

\begin{lemma}\label{hopf-surfaces} The Hopf surfaces with canonical bundle of real type are precisely those described in Theorem~\ref{main}. In particular, by Proposition~\ref{obstruction}, these are the only Hopf surfaces which can possibly admit a compatible strongly bihermitian metric.
\end{lemma}
\begin{proof}  This is rather standard (see e.g. \cite[p.669]{dloussky}). On a Hopf surface $(M,J)$, the morphism $H^1(M, \C^*) \to H^1_0(M, \cO^*)$ is injective (it is also surjective~\cite[II,p.699]{kodaira}). This is because flat $\C^*$-bundles are in bijection with $\C^*$-representations of $\Gamma=\pi_1(M)$ at one hand, and a non-vanishing holomorphic function $f$ on $\C^2$, which satisfies $f \circ \gamma_0 = \alpha f, \ \alpha \in \C^*$ (where $\gamma_0$ is the contraction \eqref{contraction}) must be constant, at the other. It follows that a flat $\C^*$-bundle is of real type if and only if the corresponding $\C^*$-representation of $\Gamma$ takes values in $\R_+^*$. The pull-back of the canonical bundle $\cK$ to the universal covering space $\C^2\setminus\{(0,0)\}$ is trivialized by the holomorphic $(2,0)$-form $\Omega= dz_1\wedge dz_2$. The contraction $\gamma_0$ acts by $\gamma_0^*(\Omega) = (\alpha \beta) \Omega$,  while for any element $h$ of the finite group $H\subset U(2)$ we have $h^*(\Omega)= {\rm det}(h)\Omega$. It follows that $\cK^{-1}$ is a flat $\C^*$-bundle associated to the representation $\gamma_0\mapsto \alpha \beta, \ h \mapsto {\rm det}(h)$; it is therefore of real type if and only if $\alpha\beta \in \R_+^*$ and $H \subset SU(2)$.

From this and Kato's classification mentioned above, the cases (a), (b) and (c) of Theorem~\ref{main} follow easily: The particular form of the contraction $\gamma_0$ is obtained from \eqref{contraction} by putting $\alpha\beta =a \in \R_+^*$. There is nothing more to prove about $H$ in the case (a) (which corresponds to $|\alpha|=|\beta|$ in Kato's classification). For the cases (b) and (c) we have $|\alpha| \neq |\beta|, \lambda =0$ and $\lambda \neq 0$, respectively. From Kato's classification $H\subset S(U(1)\times U(1))$. As $H\subset SU(2)$, it acts freely on the unit sphere $S^3 \subset \C^2$ (because the action of $SU(2)$ is the same as $S^3$ acting on itself by left multiplications). It follows that $H$ acts freely on $S^1\times \{0\}$ and $\{ 0\} \times S^1$. Thus, the projections of $U(1)\times U(1)$ to its factors inject $H\subset U(1)\times U(1)$ into $S^1$. Since $H$ is finite, it must be cyclic. There is a unique finite cyclic sub-group of $S(U(1)\times U(1))$ of order $\ell$, namely  the group generated by $(z_1,z_2)\mapsto (\varepsilon z_1, \varepsilon^{-1} z_2)$ where $\varepsilon$ is a primitive $\ell$-th root of $1$.

Finally, as the case (c) corresponds to $\lambda \neq 0$ in Kato's classification,  $\Gamma$ must be abelian, meaning that $(z_1,z_2)\mapsto (\varepsilon z_1, \varepsilon^{-1} z_2)$ (with $\varepsilon$  primitive $\ell$-th root of $1$) commutes with $\gamma_0$. This places the constraint $m = k\ell -1$,  $k\in {\mathbb N}^*$. \end{proof}

We now turn to the construction of strongly bihermitian metrics on the Hopf surfaces listed in Theorem~\ref{main}.
\begin{prop}\label{construction} Any Hopf surface $(M,J)$ described in Theorem~\ref{main} admits a strongly bihermitian metric $(g,J_+,J_-)$ with $J_+=J$.
\end{prop}
\begin{proof} The Hopf surfaces described in the case (a) of the theorem are hyperhermitian and the complex structure $J$ belongs to the underlying hypercomplex family (cf.~\cite{boyer, pontecorvo}). To see this directly, notice that in this case the fundamental group respects the standard (flat) hyperhermitian conformal structure on $\C^2$. Notice also that any hyperhermitian metric on a compact complex surface can be deformed to obtain non-hyperhermitian strongly bihermitian metrics~\cite[Prop.~1]{AGG}.

Let us now consider the Hopf surfaces described in the case (b)  of Theorem~\ref{main}.  Recall that in terms of the Kato classification referred to earlier, this is the case when $\lambda=0$ and $\alpha \beta = a \in \R_+^*$ in \eqref{contraction}, and $H$ is a finite cyclic subgroup of  $S(U(1)\times U(1))$. 

We shall use Corollary~\ref{hitchin-lemma} (with $\tau=0$) first on the universal covering space ${\C^2 \setminus \{(0,0)\}}$, then (with $\tau \neq 0$) on the quotient $(M,J)= \Big(\C^2 \setminus \{(0,0)\}\Big)/ \Gamma$.

Following \cite[Prop.~3]{hitchin4}, consider the closed $(2,0)$-form $\Omega=d{z}_1 \wedge  d {z}_2= \Phi + i \Psi$ on $\C^2\setminus\{(0,0)\}$ and a function $f$ which is a K\"ahler potential (i.e. $dd^c f>0$). Let $X$ be the $\Phi$-hamiltonian vector field corresponding to $f$,  defined by  $i_X \Phi = df$. Put $\Psi_-^t := \phi_t^* (\Psi)$, where $\phi_t$ is the flow of $X$. Then the forms $(\Phi, \Psi, \Psi_-^t)$ satisfy the conditions of Corollary~\ref{hitchin-lemma} (with $\tau=0$), except the last (positivity) condition. However,  we have $(\Psi_-^0)^{1,1} =0$ and
\begin{equation}\label{positivity}
\Big(\frac{\partial}{\partial t} \Psi_-^t\Big)_{t=0} = \cL_X \Psi= d(i_X \Psi)= -d (i_{JX} \Phi)= dd^cf,
\end{equation}
so $(\Psi_-^t)^{1,1}$ is positive definite for small $t$ (at least in a neighborhood of each point).

We now have to make a special choice of the K\"ahler potential $f$ so that the construction descends to the quotient of $\C^2 \setminus \{(0,0)\}$ by $\Gamma$. In fact, we shall use the K\"ahler potentials  introduced by Gauduchon--Ornea~\cite{GO} in order to construct (explicit) locally conformally K\"ahler metrics on Hopf surfaces. Following \cite[Remark~2]{GO}, for any $\alpha, \beta$ (as in \eqref{contraction} with $\lambda=0$) consider the flow $\varphi_t(z_1,z_2)=(\alpha^t z_1, \beta^t z_2)$ of a vector field generating the contraction, i.e.  $\gamma_0= \varphi_1$ (this requires a choice of ${\rm arg}(\alpha)$ and ${\rm arg}(\beta)$). We then define a function  $r_{\alpha,\beta} : \C^2 \setminus \{(0,0)\} \to \R$ as follows: $r_{\alpha, \beta}(z)$ is the real number such that $\varphi_{-r_{\alpha,\beta}}(z)$ belongs to the unit sphere $S^3\subset \C^2$. It is shown in \cite{GO} that this definition is correct and that the positive real function $f_{\alpha, \beta} = {\rm exp}\Big((\ln|\alpha| + \ln|\beta|) r_{\alpha,\beta} \Big)$ is a K\"ahler potential on $\C^2 \setminus \{(0,0)\}$. In fact, as $r_{\alpha,\beta}(\gamma_0 \cdot z)=  r_{\alpha,\beta}(z) +1$ by definition, we have
\begin{equation}\label{rescaling}
f_{\alpha, \beta}(\gamma_0 \cdot z)= |\alpha||\beta| f_{\alpha, \beta}(z),
\end{equation}
By the explicit formulae in \cite{GO}, $f_{\alpha,\beta}$ is also $U(1)\times U(1)$-invariant (and is $U(2)$-invariant when $|\alpha|=|\beta|$). It follows that $$F_{\alpha,\beta}:= \frac{1}{f_{\alpha,\beta}} dd^c f_{\alpha,\beta}$$ is a $\Gamma$-invariant positive definite $(1,1)$-form which defines a locally conformally K\"ahler metric on the corresponding Hopf surface $(M,J)$.

We use the function $f_{\alpha, \beta}$  in our construction of the $2$-forms $(\Phi, \Psi, \Psi_-^t)$ as above.
First of all, because in our case $\alpha \beta = |\alpha||\beta|$ and $H\subset S(U(1)\times U(1))$, it follows from \eqref{rescaling} that the corresponding vector field $X$ is $\Gamma$-invariant, and thus $X$ and its flow $\phi_t$ are defined on $(M,J)$. Secondly, by \eqref{rescaling}, the $2$-forms $(\frac{1}{f_{\alpha,\beta}}\Phi, \frac{1}{f_{\alpha,\beta}}\Psi, \frac{1}{f_{\alpha,\beta}} \Psi_-^t)$ are $\Gamma$-invariant, and therefore define a triple of $2$-forms $(\check \Phi, \check \Psi, {\check \Psi}_-^t)$ on $(M,J)$, satisfying  the conditions of Corollary~\ref{hitchin-lemma} (with $\tau = -d \log f_{\alpha,\beta}$), except possibly the last (positivity) condition. But  as $f_{\alpha,\beta}$ is invariant under $\phi_t$ by construction, we get by \eqref{positivity}  $$\Big(\frac{\partial}{\partial t}  {\check \Psi}_-^t\Big)_{t=0} = F_{\alpha,\beta},$$
which is a positive definite $(1,1)$-form on the compact manifold $(M,J)$. As the $(1,1)$-part of ${\check \Psi}_-^t$ at $t=0$ is zero, it must be positive definite for all small $t$ different than $0$.

To conclude the proof we consider the Hopf surfaces in the case (c) of Theorem~\ref{main}.  These correspond to $\lambda \neq 0$ and $\alpha=\beta^m$ in \eqref{contraction}, with $\beta^{m+1}= |\beta|^{m+1}=a$. Recall that in this case $H$ is the cyclic group generated by $(z_1,z_2)\mapsto (\varepsilon z_1, \varepsilon^{-1} z_2)$ with $\varepsilon$ primitive $\ell$-th root of $1$ so that $H$ commutes with the contraction \eqref{contraction}. 

 It is well-known  (see e.g.~\cite{harvey-lawson,KM}) that any two Hopf surfaces with fundamental groups  $\Gamma_{\beta,m,\lambda, \varepsilon}$ and $\Gamma_{\beta, m, \lambda', \varepsilon}$ corresponding to generators with the same values of $\beta$, $m$ and $\varepsilon$  (but possibly different non-zero values  $\lambda$ and $\lambda'$) are isomorphic as complex manifolds.  Because of this fact,  in order to adapt the construction in the case (b) to the case (c), it suffices to find {\it some} $\lambda_0 \neq 0$, and a {positive} function $f_{\beta,m}$ on $\C^2\setminus \{ (0,0) \}$ such that
\begin{enumerate}
\item[$\bullet$] $f_{\beta,m}(\gamma_0 \cdot z)= |\beta|^{m+1}f_{\beta,m}(z)$, where $\gamma_0$ is given by \eqref{contraction} with $\lambda=\lambda_0$;
\item[$\bullet$] $dd^c f_{\beta,m} >0$;
\item[$\bullet$] $f_{\beta,m}$ is $H$-invariant.
\end{enumerate} 
The argument in \cite[p.1125]{GO} produces positive  functions which satisfy the first two but not, a priori, the third requirement. However, their construction starts by producing a family of positive smooth functions $f_{\beta,m, \lambda}$ satisfying the first condition and  $\lim_{\lambda \to 0} f_{\beta,m, \lambda} =f_{\beta^m,\beta}$, where $f_{\beta^m,\beta}$ is the function $f_{\alpha,\beta}$ with $\alpha= \beta^m$, defined above. Since $f_{\beta^m,\beta}$ is $H$-invariant (it is, in fact, $U(1)\times U(1)$-invariant as one can see from  the equation (10) in \cite{GO}) and since $H$ commutes with $\gamma_0$, by  replacing $f_{\beta,m,\lambda}$ with its average over $H$, we can assume without loss that $f_{\beta,m,\lambda}$ are $H$-invariant. Thus, still following \cite{GO},  the $(1,1)$-form $F_{\beta,m,\lambda}:=\frac{dd^cf_{\beta,m,\lambda}}{f_{\beta,m,\lambda}}$ descends to the respective complex manifold $(M_{\lambda},J_{\lambda})$. By identifying $M_{\lambda}$ with $M_0$ (as smooth manifolds), we have $\lim_{\lambda \to 0} J_{\lambda} = J_{0}, \ \lim_{\lambda \to 0}F_{\beta,m,\lambda} = F_{\beta^m,\beta}$, where $J_{0}$ is the complex structure on $M_0$ obtained by taking $\gamma_0$ with $\alpha=\beta^m,  \lambda=0$ (note that $H$ stays unchanged). As $F_{\beta^m,\beta}$ is a  positive definite $(1,1)$-form with respect to $J_{0}$, so is $F_{\beta,m,\lambda}$ (with respect to $J_{\lambda}$) for  $\lambda$ sufficiently small. \end{proof}


\begin{thebibliography}{99}

\bibitem{A} V.~Apostolov, {\it Bihermitian surfaces with odd first Betti number}, Math. Z. {\bf 238} (2001), 555--568.

\bibitem{AG0} V.~Apostolov and  P.~Gauduchon, {\it The Riemannian Goldberg--Sachs Theorem}, Internat. J. Math. {\bf 8} (1997), 421--439.

\bibitem{AGG} V.~Apostolov, P.~Gauduchon, and G.~Grantcharov, {\it Bihermitian structures on complex surfaces}, Proc. London Math. Soc. (3) {\bf 79} (1999), 414--428. {\it Corrigendum},  {\bf 92} (2006), 200--202.

\bibitem{AG} V.~Apostolov and M.~Gualtieri, {\it Generalized K\"ahler manifolds, Commuting Complex Structures, and Split Tangent Bundles}, Comm. Math. Phys. {\bf 271} (2007), 561--575.


\bibitem{bpv} W. Barth, K.~Hulek, C. Peters, and A. Van de Ven, {\rm Compact Complex Surfaces},  Springer, Heidelberg, Second Edition, 2004.



\bibitem{bogomolov} F. A. Bogomolov, {\it Classification of surfaces of
class ${\rm VII}_0$ and affine geometry}, Math. USSR-Izv, {\bf 10} (1976), 255--269.

\bibitem{boyer} C.~P.~Boyer, {\it A note on hyperhermitian four-manifolds}, Proc. Amer. Math. Soc. {\bf 102} (1988), 157--164.



\bibitem{gualtieri-et-al}  H.~Bursztyn, G.~ R.~Cavalcanti and  M.~Gualtieri, {\it Reduction of Courant algebroids and generalized complex structures}, Adv. in Math. {\bf 211} (2007), 726--765.



\bibitem{dloussky} G.~Dloussky, {\it On surfaces of class ${\rm VII}_0^+$ with numerically anticanonical divisor}, Amer. J. Math. {\bf 128} (2006), 639--670.


\bibitem{physicists} S.~J.~Gates, C.~M.~Hull, and M.~Ro\u{c}ek, {\it Twisted multiplets and new supersymmetric nonlinear sigma models}, Nuc. Phys. B {\bf 248} (1984), 157--186.

\bibitem{gauduchon0} P.~Gauduchon, {\it Le th\'eor\`eme de dualit\'e pluriharmonique}, C.R. Acad. Sci. Paris {\bf 293} (1981), 59--63.

\bibitem{gauduchon1} P.~Gauduchon, {\it La 1-forme de torsion d'une vari\'et\'e hermitienne compacte}, Math. Ann. {\bf 267} (1984), 495--518.

\bibitem{gauduchon2} P.~Gauduchon, {\it Le premier espace de cohomologie de deRham d'une surface complexe \`a premier nombre de Betti impair}, unpublished preprint.

\bibitem{GO} P.~Gauduchon and L.~Ornea, {\it Locally conformally K\"ahler metrics on Hopf surfaces}, Ann. Inst. Fourier (Grenoble), {\bf 48} (1998), 1107--1127.

    
\bibitem{goto} R.~Goto, {\it Deformations of generalized complex and generalized K\"ahler structures}, preprint 2007, available at arxiv:math.DG/0705.2495.


\bibitem{gualtieri} M.~Gualtieri, {\it Generalized complex geometry}, D. Phil. Thesis, University of Oxford, 2003, math.DG/0401221. 

\bibitem{harvey-lawson} R.~Harvey and H. Blaine Lawson, Jr, {\it An intrinsic characterisation of K\"ahler manifolds}, Inv. Math. {\bf 74} (1983), 139--150.

\bibitem{hitchin} N.~J.~Hitchin, {\it Generalized Calabi-Yau manifolds}, Q. J. Math. {\bf 54} (2003), 281--308.

\bibitem{hitchin2} N.~J.~Hitchin, {\it Instantons and generalized K\"ahler geometry},  Comm. Math. Phys. {\bf 265} (2006), 131--164. 

\bibitem{hitchin4} N.~J.~Hitchin, {\it Bihermitian metrics on Del Pezzo Surfaces}, preprint 2006, math.DG/0608213.

\bibitem{inoue} M.~Inoue, {\it On surfaces of class ${\rm VII}_0$}, Invent. Math. {\bf 24} (1974), 269--310.

\bibitem{kato} M. Kato, {\it Topology of Hopf surfaces}, J. Math. Soc. Japan {\bf 27} (1975),  222--238. {\it Erratum} J. Math. Soc. Japan {\bf 41} (1989), 173--174.

\bibitem{kato1}M. Kato, {\it Compact Differentiable $4$-folds with Quaternionic structures}, Math. Ann. {\bf 248} (1980), 79--96. {\it Erratum}, Math. Ann. {\bf 283} (1989), 352.

\bibitem{kodaira} K. Kodaira, {\it On the structure of compact complex
analytic surfaces},  II, III, Amer. J. Math. {\bf 88} (1966), 682--721; {\bf 90} (1968), 55--83.

\bibitem{KM} K.~Kodaira and J.~Morrow, Complex manifolds. Holt, Rinehart and Winston, Inc., New York-Montreal, Que.-London, 1971.

\bibitem{kobak} P.~Kobak, {\it Explicit doubly-Hermitian metrics}, Differential Geom. Appl. {\bf 10} (1999), 179--185.





\bibitem{yau2} J.~Li, S.-T.~Yau, and F.~Zheng, {\it On projectively-flat
hermitian manifolds}, Comm. Anal. Geom. {\bf 2} (1994), 103--109.


\bibitem{Lin-Tolman} Y.~Lin and S.~Tolman, {\it Symmetries in generalized K\"{a}hler geometry}, Comm. Math. Phys. {\bf 268} (2006), 199--222.


\bibitem{pontecorvo} M.~Pontecorvo, {\it Complex structures on Riemannian four-manifolds}, Math. Ann. {\bf 309}  (1997), 159--177.

\bibitem{salamon} S.~Salamon, {\it Special structures on 4-manifolds}, Riv. Mat. Univ. Parma (4) {\bf 17} (1991), 109--123.



\bibitem{teleman0} A.~Teleman, {\it Projectively-flat surfaces and Bogomolov's theorem on class ${\rm VII}_0$ 
surfaces}, Int. J. Math. {\bf 5} (1994), 253--264.

\bibitem{teleman} A.~Teleman, {\it The pseudo-effective cone of a non-K\"ahlerian surface and applications}, Math. Ann. {\bf 335} (2006), 965--989.






\end{thebibliography}
\end{document}